\documentclass[reqno,11pt]{amsart}
\usepackage[letterpaper,hmargin=1.25in,vmargin=0.9in]{geometry}

\usepackage[T1]{fontenc}
\usepackage{lmodern}
\usepackage{microtype}

\usepackage{mathtools}
\usepackage{amssymb}
\usepackage{mathrsfs}

\usepackage[unicode,hidelinks]{hyperref}

\allowdisplaybreaks[1]
\newcommand{\R}{\mathbb{R}}
\newcommand{\Z}{\mathbb{Z}}

\DeclareMathOperator{\conv}{conv}
\DeclareMathOperator{\lin}{lin}

\DeclareMathOperator{\interior}{int}
\DeclareMathOperator{\boundary}{bd}
\DeclareMathOperator{\vol}{vol}

\DeclareMathOperator{\GL}{GL}

\newcommand{\Sym}{\mathfrak S}
\newcommand{\hk}{\widehat\kappa}
\newcommand{\latG}{\mathcal G}

\theoremstyle{plain}
\newtheorem{theorem}{Theorem}[section]
\newtheorem{proposition}[theorem]{Proposition}
\newtheorem{lemma}[theorem]{Lemma}

\theoremstyle{definition}

\newtheorem{question}[theorem]{Question}

\theoremstyle{remark}

\numberwithin{equation}{section}

\title[$\kappa$-covering minima]{On $\kappa$-covering minima}

\author{Fei Xue}
\address{School of Mathematical Sciences, Nanjing Normal University, Nanjing 210046, China}
\email{05429@njnu.edu.cn}


\date{}

\hypersetup{
  pdftitle={On kappa-covering minima},
  pdfauthor={Fei Xue},
  pdfsubject={Geometry of numbers and convex geometry},
  pdfkeywords={covering minima, lattice covering, Kuhn triangulation, successive minima}
}

\begin{document}

\begin{abstract}
We introduce $\kappa$-covering minima by requiring the permutation simplices of one independent lattice frame to fit into translates of a common dilate of a convex body. We prove a lower bound for their normalized product, with coefficient $(n+1)/2^n$ for centrally symmetric bodies and an explicit dependence on Minkowski asymmetry in general. On Euclidean balls, the supremum over full-rank lattices of the ratio of the resulting fitting radius to the covering radius has order $\sqrt n$. We also establish basic comparisons with classical lattice minima and bounds under products, projections, Minkowski addition, and free direct sums. We extend the definition by allowing projected periodic lattice triangulations. These extended minima agree with the classical covering minima on ellipsoids.
\end{abstract}

\maketitle

\section{Introduction}

A useful sufficient condition for lattice covering is to fit the simplices of a periodic lattice triangulation into translates of the covering body. Balitskiy's triangulation lemma~\cite{Balitskiy2020} gives this implication for general convex bodies. Lian and Xue \cite[Theorem~1.1]{LianXueMinimal2026} used the permutation simplices of a common lattice frame to formulate a lattice covering criterion. The same fitting radius occurs in their work on the regular octahedron~\cite[Section~9]{LianXueOctahedron2026}. Grozev and Nikolov~\cite{GrozevNikolov2026} apply the same triangulation principle to parallelotopes containing prescribed vertices of the unit cube. We define the associated minima in all dimensions and examine their relation to classical lattice minima.

Throughout, $K\subset\R^n$ is a convex body, that is, a compact convex set with nonempty interior, and $\Lambda\subset\R^n$ is a full-rank lattice. For a nonempty compact set $A\subset\R^n$, put
$$r_K(A)=\inf\{t>0:A\subset x+tK\text{ for some }x\in\R^n\}.$$
Let $B=(b_1,\ldots,b_n)$ be an ordered linearly independent tuple of lattice vectors. Its permutation simplices are
$$\Delta_\sigma(B)=\conv\left\{0,b_{\sigma(1)}, b_{\sigma(1)}+b_{\sigma(2)},\ldots,\sum_{j=1}^n b_j\right\}, \qquad \sigma\in\Sym_n,$$
where $\Sym_n$ is the permutation group. We call their family the \emph{Kuhn packet} of $B$, and define
\begin{equation}\label{eq:ordinary-kappa}
    \kappa_B(K)=\max_{\sigma\in\Sym_n}r_K(\Delta_\sigma(B)), \qquad \kappa(K,\Lambda)=\inf_{B\text{ a frame in }\Lambda}\kappa_B(K).
\end{equation}
Each simplex has its own fitting translation, while the same frame is used for the entire packet. The vectors $b_1,\ldots,b_n$ may generate a proper sublattice $\Gamma_B=\sum_j\Z b_j$ of $\Lambda$. The simplices and their $\Gamma_B$-translates form the Kuhn triangulation of $\R^n$.

For a lattice subspace $L$, meaning a subspace spanned by vectors of $\Lambda$, let $K|L^\perp$ and $\Lambda|L^\perp$ denote orthogonal projections. Write $\latG(\Lambda,k)$ for the family of $k$-dimensional lattice subspaces. The \emph{$i$th $\kappa$-covering minimum} is
\begin{equation}\label{eq:kappa-minima}
    \kappa_i(K,\Lambda)= \sup_{L\in\latG(\Lambda,n-i)} \kappa(K|L^\perp,\Lambda|L^\perp),\qquad 1\le i\le n.
\end{equation}
Thus $\kappa_n=\kappa$. We use a supremum without assuming its attainment. This projection convention is also used for covering and packing minima~\cite{HenkSchymuraXue2021,KannanLovasz1988}. In particular, writing $\mu(K,\Lambda)=\inf\{t>0:tK+\Lambda=\R^n\}$, the classical covering minima satisfy
\begin{equation}\label{eq:classical-quotient}
    \mu_i(K,\Lambda)=\sup_{L\in\latG(\Lambda,n-i)} \mu(K|L^\perp,\Lambda|L^\perp).
\end{equation}

The minima satisfy $\mu_i\le\kappa_i\le i\mu_i$, and $\kappa_i=\mu_i$ for $1\le i\le\min\{2,n\}$. These comparisons and the other basic properties are proved in Section~\ref{sec:basic}.

For the normalized product of the covering minima, Gonz\'alez Merino and Schymura~\cite{GonzalezMerinoSchymura2017} proved
$$\frac{\vol_n(K)}{\det(\Lambda)}\prod_{i=1}^n\mu_i(K,\Lambda)\ge\frac1{n!},$$
and conjectured that the coefficient can be replaced by $(n+1)/2^n$ for every convex body. In dimension two, the sharp coefficient $3/4$ was established by Schnell~\cite{Schnell1995}. Our first main result gives stronger bounds for the $\kappa$-product. Write $H_j=\sum_{\ell=1}^j1/\ell$ for the $j$th harmonic number.

\begin{theorem}\label{thm:main-product}
For every convex body $K\subset\R^n$ and full-rank lattice $\Lambda$,
\begin{equation}\label{volprod:harmonic}
    \frac{\vol_n(K)}{\det(\Lambda)}\prod_{i=1}^n\kappa_i(K,\Lambda) \ge\frac1{n!}\prod_{j=1}^n H_j.
\end{equation}
If $K$ is centrally symmetric, then
\begin{equation}\label{volprod:symmetric}
    \frac{\vol_n(K)}{\det(\Lambda)}\prod_{i=1}^n\kappa_i(K,\Lambda) \ge\frac{n+1}{2^n}.
\end{equation}
\end{theorem}

Since $\mu_i\le\kappa_i$, the classical bound $1/n!$ also holds for the normalized $\kappa$-product. The bound~\eqref{volprod:harmonic} improves this coefficient by the factor $\prod_{j=1}^nH_j$. For centrally symmetric bodies, \eqref{volprod:symmetric} gives the coefficient conjectured for the classical covering product. Section~\ref{sec:volume-product} gives a refinement in terms of Minkowski asymmetry.

The common-frame requirement also has a quantitative cost. Let $B_2^n$ denote the Euclidean unit ball.

\begin{theorem}\label{thm:euclidean-separation}
There are absolute constants $c,C>0$ such that every full-rank lattice $\Lambda\subset\R^n$ satisfies
$$\kappa(B_2^n,\Lambda)\le C\sqrt n\,\mu(B_2^n,\Lambda),$$
whereas, for every sufficiently large $n$, some determinant-one lattice $\Lambda_n$ satisfies
$$\kappa(B_2^n,\Lambda_n)\ge c\sqrt n\,\mu(B_2^n,\Lambda_n).$$
Thus $\sup_\Lambda\kappa(B_2^n,\Lambda)/\mu(B_2^n,\Lambda) =\Theta(\sqrt n)$.
\end{theorem}

This shows that $\kappa$ cannot be bounded by a fixed power of $\log(2n)$ times the Kannan--Lov\'asz projection-volume parameter, which compares projected lattice determinants with the corresponding projected body volumes. Section~\ref{sec:separation} proves the theorem and this consequence.

We also establish formulas and bounds for Cartesian products, projections, Minkowski addition, and free direct sums in Section~\ref{sec:operations}. The final section enlarges the definition by allowing projected periodic lattice triangulations. The resulting minima $\hk_i$ lie between $\mu_i$ and $\kappa_i$ and equal $\mu_i$ on ellipsoids. For general convex bodies, the corresponding projection-volume bounds remain open.

\section{Basic properties}\label{sec:basic}

A body $C$ is called $o$-symmetric if $C=-C$. For such a body, its $j$th successive minimum is
$$\lambda_j(C,\Lambda)=\inf\{t>0:\dim\lin(tC\cap\Lambda)\ge j\}.$$
Write $h_K(u)=\max_{x\in K}\langle u,x\rangle$ for the support function of $K$. The dual lattice is $\Lambda^*=\{u:\langle u,\Lambda\rangle\subset\Z\}$. A surjective linear map $\pi$ is called a \emph{lattice quotient} if $\pi\Lambda$ is discrete. Its kernel $L$ is a lattice subspace, and $\pi$ is the orthogonal projection onto $L^\perp$ followed by a linear isomorphism. Since $\kappa(AK,A\Lambda)=\kappa(K,\Lambda)$ for every invertible $A$, we may use these maps in the quotient formula~\eqref{eq:kappa-minima}.

\begin{proposition}\label{prop:basic}
For $A\in\GL(n,\R)$, $x\in\R^n$, and $a,b>0$,
$$\kappa_i(AK+x,A\Lambda)=\kappa_i(K,\Lambda),\qquad \kappa_i(aK,b\Lambda)=\frac ba\kappa_i(K,\Lambda).$$
For convex bodies $K\subset H$ and full-rank lattices $\Lambda\subset\Gamma$,
$$\kappa_i(H,\Lambda)\le\kappa_i(K,\Lambda),\qquad \kappa_i(K,\Gamma)\le\kappa_i(K,\Lambda).$$
The sequence $\kappa_1,\ldots,\kappa_n$ is nondecreasing. For every rank-$r$ lattice quotient $\pi$,
$$\kappa_j(\pi K,\pi\Lambda)\le\kappa_j(K,\Lambda)\qquad(1\le j\le r).$$
\end{proposition}
\begin{proof}
For every frame $B$, a change of coordinates gives
$$\kappa_{AB}(AK+x)=\kappa_B(K),\qquad \kappa_{bB}(aK)=\frac ba\kappa_B(K).$$
Taking infima over frames proves the first two identities for $\kappa$. Enlarging the body decreases the fitting radii, while enlarging the lattice adds admissible frames. These facts pass to every index by the quotient formula: invertible maps transport lattice subspaces, and finite-index lattice inclusions have the same lattice subspaces.

For a rank-$r$ lattice quotient $\pi$ and a frame $B$, choose $I\subseteq\{1,\ldots,n\}$ with $|I|=r$ and independent images $\pi b_j$, $j\in I$. Place the indices of $I$ first in a permutation of all $n$ indices. The face on its first $r+1$ vertices projects to the corresponding permutation simplex of $\pi B_I$. Thus
$$\kappa(\pi K,\pi\Lambda)\le\kappa_{\pi B_I}(\pi K)\le\kappa_B(K).$$
Taking the infimum over $B$ proves the assertion for the ordinary radius. Composing lattice quotients proves the indexed assertion. Finally, every rank-$i$ lattice quotient factors through one of rank $i+1$, which gives $\kappa_i\le\kappa_{i+1}$.
\end{proof}

\begin{proposition}\label{prop:regularity}
There is an independent frame $B$ in $\Lambda$ such that $\kappa_B(K)=\kappa(K,\Lambda)$. Each $\kappa_i$ is jointly continuous for the Hausdorff topology on convex bodies and the usual topology on full-rank lattices. More precisely, if $K$ contains a translate of $rB_2^n$ and $\delta=d_H(K,H)<r$, then
\begin{equation}\label{eq:body-continuity}
    \frac{\kappa_i(K,\Lambda)}{1+\delta/r}\le\kappa_i(H,\Lambda)\le\frac{\kappa_i(K,\Lambda)}{1-\delta/r}.
\end{equation}
\end{proposition}
\begin{proof}
Fix a frame $B_0$ in $\Lambda$ and put $R=\kappa_{B_0}(K)$. Compactness gives attainment of $r_K(\Delta_\sigma(B))$ for each frame $B$ and permutation $\sigma$. Since every $b_j$ is an edge vector of a permutation simplex,
$$\kappa_B(K)\le R\quad\Longrightarrow\quad b_j\in\Lambda\cap R(K-K)\qquad(1\le j\le n).$$
Discreteness of $\Lambda$ leaves only finitely many such frames, including $B_0$. The minimum of $\kappa_B(K)$ over this finite family equals $\kappa(K,\Lambda)$.

Translate both bodies so that $rB_2^n\subset K$. The support-function inequalities $|h_H-h_K|\le\delta$ and $h_K\ge r$ give
$$(1-\delta/r)K\subset H\subset(1+\delta/r)K.$$
Proposition~\ref{prop:basic} gives~\eqref{eq:body-continuity}, uniformly over all lattice quotients. For $\Lambda_m\to\Lambda$, choose basis matrices $A_m\to A$. The fixed-lattice estimate proves joint continuity, since $K_m\to K$ implies
\begin{equation*}
    A_m^{-1}K_m\to A^{-1}K,\qquad \kappa_i(K_m,\Lambda_m)=\kappa_i(A_m^{-1}K_m,\Z^n).\qedhere
\end{equation*}
\end{proof}

\begin{proposition}\label{prop:successive}
For $1\le i\le n$,
\begin{align*}
    \kappa_i(K,\Lambda)&\le\sum_{j=n-i+1}^n\lambda_j(K-K,\Lambda),\\
    \kappa_i(C,\Lambda)&\le\frac12\sum_{j=n-i+1}^n\lambda_j(C,\Lambda)\qquad(C=-C).
\end{align*}
\end{proposition}
\begin{proof}
Put $\alpha_j=\lambda_j(K-K,\Lambda)$ and choose independent $b_j=\alpha_j(x_j-y_j)\in\Lambda$ with $x_j,y_j\in K$, as in the standard segment construction~\cite{KannanLovasz1988}. For every index set $I$,
$$\sum_{j\in I}[0,b_j]+\sum_{j\in I}\alpha_jy_j=\sum_{j\in I}\alpha_j[y_j,x_j]\subset\left(\sum_{j\in I}\alpha_j\right)K.$$
For a rank-$i$ lattice quotient $\pi$, choose $I$ with $|I|=i$ and independent images $\pi b_j$, $j\in I$. The projected parallelotope contains all permutation simplices of this image frame. Consequently,
$$\kappa(\pi K,\pi\Lambda)\le\sum_{j\in I}\alpha_j\le\sum_{j=n-i+1}^n\alpha_j.$$
Take the supremum over $\pi$. The identity $C-C=2C$ gives the second bound.
\end{proof}

\begin{proposition}\label{prop:comparison}
For every convex body $K\subset\R^n$ and full-rank lattice $\Lambda$,
\begin{equation}\label{eq:basic-comparison}
    0<\kappa_1(K,\Lambda)\le\cdots\le\kappa_n(K,\Lambda)<\infty,\qquad \mu_i(K,\Lambda)\le\kappa_i(K,\Lambda)\le i\mu_i(K,\Lambda).
\end{equation}
\end{proposition}
\begin{proof}
Choose an optimal frame $B$ and put $t=\kappa(K,\Lambda)$. Every simplex in its Kuhn triangulation fits into a translate of $tK$. The triangulation has uniformly bounded simplex diameters and vertex set $\Gamma_B\subset\Lambda$. Balitskiy's triangulation lemma~\cite{Balitskiy2020} gives
$$tK+\Gamma_B=\R^n,\qquad tK+\Lambda=\R^n.$$
Thus $\mu(K,\Lambda)\le\kappa(K,\Lambda)$.

For the reverse bound, translate $K$ so that $0\in K$ and set $t=\mu(K,\Lambda)$, so $tK+\Lambda=\R^n$. Fix $\lambda\in\Lambda$. Only finitely many of these translates meet the segment $[0,\lambda]$. Their intersections with the segment are closed intervals covering it. We can therefore join $tK$ to $tK+\lambda$ by a finite chain of translates whose consecutive members intersect. Each successive translation difference belongs to $t(K-K)\cap\Lambda$. Adding the differences along the chain gives $\lambda$, hence
$$\Lambda=\sum_{v\in t(K-K)\cap\Lambda}\Z v,\qquad \lambda_n(K-K,\Lambda)\le t.$$
Proposition~\ref{prop:successive} now gives
$$\kappa(K,\Lambda)\le\sum_{j=1}^n\lambda_j(K-K,\Lambda)\le n\mu(K,\Lambda).$$
Apply both bounds in every rank-$i$ lattice quotient. Positivity and finiteness follow from the classical covering minima, and index monotonicity was proved in Proposition~\ref{prop:basic}.
\end{proof}

\begin{proposition}\label{prop:first-two}
One has
$$\kappa_1(K,\Lambda)=\mu_1(K,\Lambda)=\frac1{\displaystyle\min_{0\ne u\in\Lambda^*}\bigl(h_K(u)+h_K(-u)\bigr)}.$$
If $n\ge2$, then $\kappa_2(K,\Lambda)=\mu_2(K,\Lambda)$.
\end{proposition}
\begin{proof}
An interval of length $\ell$ with lattice spacing $a$ has both radii equal to $a/\ell$. For primitive $u\in\Lambda^*$, the map $x\mapsto\langle u,x\rangle$ sends $\Lambda$ onto $\Z$ and $K$ onto an interval of length $h_K(u)+h_K(-u)$. Taking the supremum over these rank-one quotients gives the width formula.

In dimension two, the two permutation triangles of a frame are translates of a triangle and its reflection. We show that the centrally symmetric tile supplied by F\'ary's theorem contains such a pair. Put $t=\mu(K,\Lambda)$. Since $tK+\Lambda$ covers the plane, F\'ary's theorem~\cite[p.~158]{Fary1950} gives a centrally symmetric hexagon or parallelogram $H\subset tK$ tiling by $\Lambda$.

Let $p$ be a vertex of $H$. The tiles $H+\lambda$ incident to $p$ cover a neighborhood of $p$, and the points $p-\lambda$ lie on $\boundary H$. At least three of these points are noncollinear. Suppose, to the contrary, that all these points lie on a line $\ell$. If $\ell$ supports $H$, the incident tiles leave a half-plane at $p$ uncovered. Otherwise, $\ell$ crosses $\interior H$ and meets $\boundary H$ in only two points. Thus at most two tiles meet at $p$, and their total angle is less than $2\pi$, since $p$ is a vertex of at least one of them. Both cases contradict the covering of a neighborhood of $p$.

Choose three noncollinear incident translations $\lambda_0,\lambda_1,\lambda_2$ and put $b_1=\lambda_0-\lambda_1$, $b_2=\lambda_1-\lambda_2$. The triangle on $p-\lambda_0,p-\lambda_1,p-\lambda_2$ is a translate of
$$T=\conv\{0,b_1,b_1+b_2\}.$$
Its reflection in the center of $H$ is a translate of the other permutation triangle, since
$$\conv\{0,b_2,b_1+b_2\}=(b_1+b_2)-T.$$
Both fit in $H\subset tK$, so $\kappa(K,\Lambda)\le t=\mu(K,\Lambda)$. Proposition~\ref{prop:comparison} gives equality. Applying the planar identity in every rank-two quotient proves the assertion for $\kappa_2$.
\end{proof}

\section{Products, projections, and addition}\label{sec:operations}

Cartesian products admit an exact formula even for projections that mix the factors. Minkowski addition gives a harmonic inequality for a fixed frame, while free direct sums give additive bounds.

\subsection{Cartesian products}

\begin{proposition}\label{prop:product}
Let $K\subset\R^n$ and $H\subset\R^m$ be convex bodies, and let $\Lambda$ and $\Gamma$ be full-rank lattices in the respective spaces. For $1\le i\le n+m$,
\begin{equation}\label{eq:cartesian-product}
    \kappa_i(K\times H,\Lambda\oplus\Gamma) =\max\{\kappa_{\min(i,n)}(K,\Lambda), \kappa_{\min(i,m)}(H,\Gamma)\}.
\end{equation}
\end{proposition}
\begin{proof}
We first prove the formula for the ordinary radius. Translate $K$ and $H$ to contain the origin, choose optimal frames $B=(b_1,\ldots,b_n)$ and $D=(d_1,\ldots,d_m)$, and put $t=\max\{\kappa(K,\Lambda),\kappa(H,\Gamma)\}$. Use the frame
$$((b_1,0),\ldots,(b_n,0),(0,d_1),\ldots,(0,d_m))$$
in the product space. In any ordering of these vectors, deleting the steps from the second factor leaves an ordering of $b_1,\ldots,b_n$. Thus the first projection of the resulting simplex is a permutation simplex of $B$; the second projection is one of $D$. If these projections lie in $x+tK$ and $y+tH$, the original simplex lies in $(x,y)+t(K\times H)$. The two factor projections give the opposite inequality by Proposition~\ref{prop:basic}. Hence
$$\kappa(K\times H,\Lambda\oplus\Gamma)=\max\{\kappa(K,\Lambda),\kappa(H,\Gamma)\}.$$

To prove the indexed upper bound, take any rank-$i$ lattice quotient $\pi$ of the product and restrict it to the two factors:
$$\pi_1x=\pi(x,0),\qquad \pi_2y=\pi(0,y).$$
Put $U_j=\operatorname{im}\pi_j$, $C_1=\pi_1K$, $C_2=\pi_2H$, and $\Lambda_1=\pi_1\Lambda$, $\Lambda_2=\pi_2\Gamma$. Then
$$\pi(K\times H)=C_1+C_2,\qquad \pi(\Lambda\oplus\Gamma)=\Lambda_1+\Lambda_2.$$
Each $\Lambda_j$ is a full-rank lattice in $U_j$: it spans that space and is a subgroup of the discrete lattice $\pi(\Lambda\oplus\Gamma)$. If $U_2=\{0\}$, the projected radius is $\kappa(\pi_1K,\pi_1\Lambda)\le\kappa_i(K,\Lambda)$; the case $U_1=\{0\}$ is analogous. Otherwise, writing $r_j=\dim U_j$, we have
\begin{align*}
    \kappa(C_1,\Lambda_1)&\le\kappa_{r_1}(K,\Lambda)\le\kappa_{\min(i,n)}(K,\Lambda),\\
    \kappa(C_2,\Lambda_2)&\le\kappa_{r_2}(H,\Gamma)\le\kappa_{\min(i,m)}(H,\Gamma).
\end{align*}
Here $r_1\le\min(i,n)$ and $r_2\le\min(i,m)$. The image spaces may intersect. The addition map
$$U_1\times U_2\longrightarrow U_1+U_2,\qquad (u_1,u_2)\longmapsto u_1+u_2$$
is a lattice quotient because its lattice image is $\Lambda_1+\Lambda_2$. Applying quotient monotonicity and then the ordinary product formula gives
\begin{align*}
    \kappa(\pi(K\times H),\pi(\Lambda\oplus\Gamma))&=\kappa(C_1+C_2,\Lambda_1+\Lambda_2)\\
    &\le\kappa(C_1\times C_2,\Lambda_1\oplus\Lambda_2)\\
    &=\max\{\kappa(C_1,\Lambda_1),\kappa(C_2,\Lambda_2)\}.
\end{align*}
The preceding factor bounds, followed by the supremum over $\pi$, prove the upper bound in~\eqref{eq:cartesian-product}.

For the lower bound, let $a=\min(i,n)$. The projection onto the first factor and index monotonicity give
$$\kappa_a(K,\Lambda)\le\kappa_a(K\times H,\Lambda\oplus\Gamma)\le\kappa_i(K\times H,\Lambda\oplus\Gamma).$$
The second factor gives the other term in the maximum.
\end{proof}

To compute the standard examples below, we use a bound by finitely many coordinate sections. The analogous estimate for covering minima is due to Krivoku\'ca~\cite[Theorem~1.3]{Krivokuca2026}.

\begin{lemma}\label{prop:coordinate-sections}
Let $v_1,\ldots,v_n\in\Lambda$ be independent, and put $L_I=\lin\{v_j:j\in I\}$. Fix $1\le i\le n$ and assume that $K\cap L_I$ has dimension $i$ whenever $|I|=i$. Then
$$\kappa_i(K,\Lambda) \le\max_{|I|=i}\kappa(K\cap L_I,\Lambda\cap L_I).$$
\end{lemma}
\begin{proof}
For a rank-$i$ lattice quotient $\pi$, choose $I$ such that $\pi v_j$, $j\in I$, are independent. Then $\pi|_{L_I}$ is an isomorphism and
$$\pi(K\cap L_I)\subset\pi K,\qquad \pi(\Lambda\cap L_I)\subset\pi\Lambda.$$
Both lattices have full rank, so Proposition~\ref{prop:basic} gives
$$\kappa(\pi K,\pi\Lambda)\le\kappa(\pi(K\cap L_I),\pi(\Lambda\cap L_I))=\kappa(K\cap L_I,\Lambda\cap L_I).$$
Take the supremum over $\pi$.
\end{proof}

The following examples show that $\kappa_i$ and $\mu_i$ agree on the standard balls and simplex for the coordinate lattice. Section~\ref{sec:separation} shows how this changes for Euclidean balls with other lattices. Let $B_p^n=\{x:\sum_j|x_j|^p\le1\}$ for $1\le p<\infty$, let $B_\infty^n=[-1,1]^n$, and put $T_n=\conv\{0,e_1,\ldots,e_n\}$. For $1\le i\le n$,
\begin{equation}\label{eq:standard-examples}
    \begin{split}
        \kappa_i(B_p^n,\Z^n)&=\mu_i(B_p^n,\Z^n)=\frac{i^{1/p}}2 \qquad(1\le p\le\infty),\\
        \kappa_i(T_n,\Z^n)&=\mu_i(T_n,\Z^n)=i,
    \end{split}
\end{equation}
where $i^{1/\infty}=1$. To see the upper bounds, apply Lemma~\ref{prop:coordinate-sections} to the standard basis. Each $i$-dimensional coordinate section is $B_p^i$ or $T_i$. The coordinate Kuhn simplices lie in $[0,1]^i$, and
$$[-1/2,1/2]^i\subset(i^{1/p}/2)B_p^i,\qquad [0,1]^i\subset iT_i.$$
For the reverse bounds, project onto $i$ coordinate directions. Nearest-integer distances give $\mu(B_p^i,\Z^i)=i^{1/p}/2$, and the supremum of the sum of fractional parts gives $\mu(T_i,\Z^i)=i$. Now use $\mu_i\le\kappa_i$.

\subsection{Minkowski addition}

For convex bodies $K,H\subset\R^n$, the sum $K+H$ contains a translate of each factor. Hence
$$\kappa_i(K+H,\Lambda) \le\min\{\kappa_i(K,\Lambda),\kappa_i(H,\Lambda)\} \qquad(1\le i\le n).$$
A stronger harmonic inequality holds when the frame is fixed.

\begin{lemma}\label{prop:fixed-frame-harmonic}
For every independent frame $B$ and convex bodies $K,H$,
$$\frac1{\kappa_B(K+H)} \ge\frac1{\kappa_B(K)}+\frac1{\kappa_B(H)}.$$
\end{lemma}
\begin{proof}
Put $r=\kappa_B(K)$ and $s=\kappa_B(H)$. The fitting radii are attained, so for each $\Delta=\Delta_\sigma(B)$ there are $x,y$ with $x+\Delta\subset rK$ and $y+\Delta\subset sH$. For every $z\in\Delta$, adding $(x+z)/r\in K$ and $(y+z)/s\in H$ gives
$$x/r+y/s+(1/r+1/s)\Delta\subset K+H.$$
Thus $\kappa_B(K+H)\le(1/r+1/s)^{-1}$, which proves the claim.
\end{proof}

The optimal frames for $K$ and $H$ can differ, so the fixed-frame inequality does not imply the same inequality for $\kappa$. In the first two indices, equality with the covering minima gives the following sharp constants.

\begin{proposition}\label{prop:low-index-addition}
For $K,H\subset\R^n$ and $1\le i\le\min\{2,n\}$,
$$\frac1{\kappa_i(K+H,\Lambda)} \ge c_i\left(\frac1{\kappa_i(K,\Lambda)} +\frac1{\kappa_i(H,\Lambda)}\right), \qquad c_1=1,\quad c_2=\frac34.$$
Neither constant can be increased uniformly over all dimensions.
\end{proposition}
\begin{proof}
For $i=1$, use Proposition~\ref{prop:first-two} and the additivity of directional widths under Minkowski addition. The minimum of their sum over nonzero dual lattice vectors is at least the sum of their separate minima.

For $i=2$, Codenotti, Freyer, and Krivoku\'ca \cite[Theorem~3.4]{CodenottiFreyerKrivokuca2026} prove
$$\frac1{\mu(A+B,\Gamma)} \ge\frac34\left(\frac1{\mu(A,\Gamma)} +\frac1{\mu(B,\Gamma)}\right)$$
for arbitrary planar convex bodies and a common full-rank lattice. For any convex body $X\subset\R^n$, Proposition~\ref{prop:first-two} and the quotient formula give
$$\frac1{\kappa_2(X,\Lambda)}=\inf_\pi\frac1{\mu(\pi X,\pi\Lambda)},$$
where $\pi$ ranges over all rank-two lattice quotients. Apply the planar inequality to $A=\pi K$ and $B=\pi H$ for the same $\pi$. The two terms on its right are at least $1/\kappa_2(K,\Lambda)$ and $1/\kappa_2(H,\Lambda)$, respectively. Taking the infimum of its left-hand side over $\pi$ proves the assertion.

Intervals give equality for $c_1$. For $c_2$, the example from \cite[equation~(1.2)]{CodenottiFreyerKrivokuca2026} is
$$P_1=[0,(1,1)]+[0,(-1/2,1/2)],\qquad P_2=[0,(-1,1)]+[0,(1/2,1/2)].$$
These parallelograms cover with $\Z^2$ and have area one, so their covering radii, and hence their packet radii, are one. Their sum is a translate of $(3/2)B_1^2$, whose radius is $2/3$ by \eqref{eq:standard-examples}. Thus equality holds for $c_2=3/4$. The coefficient cannot be replaced by one, because
$$\frac1{\kappa(P_1+P_2,\Z^2)}=\frac32<2=\frac1{\kappa(P_1,\Z^2)}+\frac1{\kappa(P_2,\Z^2)}.$$
\end{proof}

\subsection{Free direct sums}

For convex bodies $K\subset\R^n$ and $H\subset\R^m$ containing the origin, their \emph{free direct sum} is
$$K\oplus H=\conv\bigl((K\times\{0\}) \cup(\{0\}\times H)\bigr).$$
Here $\Lambda\oplus\Gamma$ is the direct sum of the lattices. The factors may contain the origin on their boundary. For the covering radius, Codenotti, Santos, and Schymura~\cite[Corollary~2.2]{CodenottiSantosSchymura2022} proved
$$\mu(K\oplus H,\Lambda\oplus\Gamma)=\mu(K,\Lambda)+\mu(H,\Gamma).$$
We compare $\kappa$ under the same operation.

\begin{proposition}\label{prop:free-sum-sandwich}
For such $K,H$ and full-rank lattices $\Lambda,\Gamma$,
\begin{align}
    &\max\{\kappa(K,\Lambda),\kappa(H,\Gamma), \mu(K,\Lambda)+\mu(H,\Gamma)\} \notag\\[-2pt]
    &\hspace{15mm}\le\kappa(K\oplus H,\Lambda\oplus\Gamma) \le\kappa(K,\Lambda)+\kappa(H,\Gamma). \label{eq:free-sum-sandwich}
\end{align}
If each factor satisfies $\kappa=\mu$, then
$$\kappa(K\oplus H,\Lambda\oplus\Gamma) =\kappa(K,\Lambda)+\kappa(H,\Gamma).$$
\end{proposition}
\begin{proof}
Put $r=\kappa(K,\Lambda)$ and $s=\kappa(H,\Gamma)$, and choose optimal factor frames. Embed them in the two coordinate subspaces, as in the proof of Proposition~\ref{prop:product}. Each permutation simplex of the resulting frame fits in a translate of $rK\times sH$. For $x\in K$ and $y\in H$,
$$\frac{(rx,sy)}{r+s}=\frac r{r+s}(x,0)+\frac s{r+s}(0,y)\in K\oplus H.$$
Hence $rK\times sH\subset(r+s)(K\oplus H)$, proving the upper bound. Factor projections give the first two lower bounds. The displayed covering-radius identity and $\mu\le\kappa$ give the third. If $\kappa=\mu$ in both factors, the upper and lower bounds coincide.
\end{proof}

\section{Products of the minima}\label{sec:volume-product}

Since $\mu_i\le\kappa_i$, a lower bound for the $\kappa$-product does not imply the conjectured bound for the covering product. We refine Theorem~\ref{thm:main-product} using the Minkowski asymmetry
$$s(K)=\inf\{s\ge1:-(K-z)\subset s(K-z)\text{ for some }z\in\interior K\}.$$
The infimum is attained. Indeed, for $s_m\downarrow s(K)$, the corresponding centers have a convergent subsequence in $K$, and the limiting center satisfies the inclusion. It lies in $\interior K$, since a supporting hyperplane at a boundary center would put the body and its reflection on opposite sides.

\begin{theorem}\label{volprod:theorem}
For every convex body $K\subset\R^n$ and full-rank lattice $\Lambda$,
\begin{equation}\label{volprod:asymmetry}
    \frac{\vol_n(K)}{\det(\Lambda)}\prod_{i=1}^n\kappa_i(K,\Lambda)\ge c_n(s(K)),\qquad c_n(s)=\prod_{j=1}^n\max\left\{\frac{H_j}{j},\frac{j+1}{j(1+s)}\right\}.
\end{equation}
In particular, Theorem~\ref{thm:main-product} holds.
\end{theorem}
The proof uses two sums of lattice segments. Averaging the permutation simplices by the Hunt--Dyson identity~\cite{Gioev2002} gives the harmonic term, while taking differences of their vertices gives the term involving $s(K)$. Both constructions are estimated using the following projection lemma.

\begin{lemma}\label{lem:zonotope-projection}
Let $K\subset\R^d$, $d\ge2$, be a convex body and $\Lambda$ a full-rank lattice. Suppose a translate of $Z=\sum_{j=1}^m\alpha_j[0,g_j]$ lies in $K$, where $\alpha_j\ge0$ and $g_j\in\Lambda\setminus\{0\}$. Then
\begin{equation}\label{volprod:segment}
    \frac{\vol_d(K)}{\det(\Lambda)}\ge\frac1d\sum_{j=1}^m\alpha_j\frac{\vol_{d-1}(K|g_j^\perp)}{\det(\Lambda|g_j^\perp)}.
\end{equation}
\end{lemma}
\begin{proof}
For a convex body $D$, integration along the fibers parallel to $g$ gives
$$\vol_d(D+u[0,g])=\vol_d(D)+u\|g\|_2\vol_{d-1}(D|g^\perp)\qquad(u\ge0).$$
Adding the segments of $\varepsilon Z$ one at a time and letting $\varepsilon\downarrow0$, we obtain
$$\lim_{\varepsilon\downarrow0}\frac{\vol_d(K+\varepsilon Z)-\vol_d(K)}{\varepsilon}=\sum_{j=1}^m\alpha_j\|g_j\|_2\vol_{d-1}(K|g_j^\perp).$$
If $x+Z\subset K$, then $K+\varepsilon Z+\varepsilon x\subset(1+\varepsilon)K$, so this limit is at most $d\vol_d(K)$. For every nonzero lattice vector $g$,
$$\|g\|_2\det(\Lambda|g^\perp)\ge\det(\Lambda),$$
as follows by extending the primitive lattice vector in its direction to a basis. Dividing by $d\det(\Lambda)$ proves~\eqref{volprod:segment}.
\end{proof}

\begin{proof}[Proof of Theorem~\ref{volprod:theorem}]
Translate $K$ so that an optimal center is at the origin and put $r=s(K)$, so $-K\subset rK$. Write
$$a_d(r)=\max\left\{\frac{H_d}{d},\frac{d+1}{d(1+r)}\right\},\qquad c_n(r)=\prod_{d=1}^n a_d(r).$$
We obtain the two coefficients from two sums of lattice segments, then apply the resulting projection estimate successively. For $n\ge2$, choose an optimal frame $B=(b_1,\ldots,b_n)$ and put $t=\kappa(K,\Lambda)$ and $C=tK$. Set $b_\varnothing=0$ and
$$b_I=\sum_{j\in I}b_j,\qquad b_*=\sum_{j=1}^n b_j,\qquad w_I=\frac1{|I|\binom n{|I|}}\quad(\varnothing\ne I\subseteq\{1,\ldots,n\}).$$

First, applying the classical Hunt--Dyson identity~\cite{Gioev2002} to the scalars $\langle u,b_j\rangle$ shows that the support functions of both sides below agree:
\begin{equation}\label{volprod:average}
    \frac1{n!}\sum_{\sigma\in\Sym_n}\Delta_\sigma(B)=\sum_{I\ne\varnothing}w_I[0,b_I].
\end{equation}
Choose $x_\sigma$ with $\Delta_\sigma(B)\subset x_\sigma+C$. By convexity,
$$\frac1{n!}\sum_\sigma\Delta_\sigma(B)\subset\frac1{n!}\sum_\sigma x_\sigma+C.$$
Thus the sum of segments in~\eqref{volprod:average} has a translate in $C$. All $b_I$ in the sum are nonzero, and $\sum_Iw_I=H_n$. Lemma~\ref{lem:zonotope-projection} therefore gives at least one lattice direction whose normalized projection volume is at most $n/H_n$ times $\vol_n(C)/\det(\Lambda)$.

Second, every $b_I$ is the difference of two vertices of a fitted permutation simplex, so $\pm b_I\in C-C$. The vertices of $P_B=\sum_j[0,b_j]$ are the subset sums $b_I$. Since $P_B-b_*=-P_B$,
$$P_B+[0,b_*]-b_*=\conv(P_B\cup(-P_B))\subset C-C\subset(1+r)C.$$
Thus a translate of $(P_B+[0,b_*])/(1+r)$ lies in $C$. This is a sum of $n+1$ nonzero lattice segments with total weight $(n+1)/(1+r)$. Lemma~\ref{lem:zonotope-projection} gives the corresponding projection-volume coefficient $(n+1)/(n(1+r))$.

Taking the construction with the larger coefficient yields a rank-$(n-1)$ lattice projection $\pi$ with
$$\frac{\vol_n(C)}{\det(\Lambda)}\ge a_n(r)\frac{\vol_{n-1}(\pi C)}{\det(\pi\Lambda)}.$$
Since $C=tK$, this is
\begin{equation}\label{volprod:one-step}
    \kappa(K,\Lambda)\frac{\vol_n(K)}{\det(\Lambda)}\ge a_n(r)\frac{\vol_{n-1}(\pi K)}{\det(\pi\Lambda)}.
\end{equation}

Starting with $K_n=K$ and $\Lambda_n=\Lambda$, apply~\eqref{volprod:one-step} successively to obtain $K_{d-1}=\pi_dK_d$ and $\Lambda_{d-1}=\pi_d\Lambda_d$ for $d=n,\ldots,2$. Each projected body still satisfies $-K_d\subset rK_d$. In dimension one,
$$\kappa(K_1,\Lambda_1)\frac{\vol_1(K_1)}{\det(\Lambda_1)}=1=a_1(r).$$
Multiplying the estimates cancels all intermediate normalized volumes. Moreover, each $K_d$ is a rank-$d$ lattice quotient of the original body, so $\kappa(K_d,\Lambda_d)\le\kappa_d(K,\Lambda)$. Consequently,
\begin{align*}
    \frac{\vol_n(K)}{\det(\Lambda)}\prod_{d=1}^n\kappa_d(K,\Lambda)&\ge\frac{\vol_n(K)}{\det(\Lambda)}\prod_{d=1}^n\kappa(K_d,\Lambda_d)\\
    &\ge\prod_{d=1}^n a_d(r)=c_n(r).
\end{align*}
Since $r=s(K)$, this proves~\eqref{volprod:asymmetry}; the interval identity also covers $n=1$. The harmonic terms give~\eqref{volprod:harmonic}. For centrally symmetric bodies, $s(K)=1$ and $H_d\le(d+1)/2$ give $c_n(1)=(n+1)/2^n$, proving~\eqref{volprod:symmetric}.
\end{proof}

\section{Euclidean separation}\label{sec:separation}

We now prove Theorem~\ref{thm:euclidean-separation}, which quantifies the difference between fitting one common frame and covering by lattice translates. The lower bound comes from the subset sums of every frame. For the upper bound, we construct a frame from the facet vectors of the Voronoi cell.

\begin{proof}[Proof of Theorem~\ref{thm:euclidean-separation}]
Write $\lambda_1(\Lambda)$ for the length of a shortest nonzero vector. Choose an optimal frame $B=(b_1,\ldots,b_n)$ and put $t=\kappa(B_2^n,\Lambda)$. Each subset sum $s_A=\sum_{j\in A}b_j$ joins two vertices of one fitted simplex, obtained by placing the indices of $A$ first. Hence $\|s_A\|_2\le2t$. Averaging over all subsets gives
$$4t^2\ge\mathbb E\|s_A\|_2^2 =\frac14\left(\sum_j\|b_j\|_2^2+ \Bigl\|\sum_jb_j\Bigr\|_2^2\right) \ge\frac n4\lambda_1(\Lambda)^2.$$
Consequently,
\begin{equation}\label{eq:short-euclidean-lower}
    \kappa(B_2^n,\Lambda)\ge\frac{\sqrt n}{4}\lambda_1(\Lambda).
\end{equation}

By Shapira and Weiss~\cite[Proposition~3]{ShapiraWeiss2014}, for an absolute $C_1>0$ and all sufficiently large $n$, there is a determinant-one lattice $\Lambda_n$ such that
$$h_n=\left(\frac n{C_1}\right)^{(n-1)/(2n)}, \qquad \lambda_1(\Lambda_n)\ge h_n,$$
and
$$\det(\Delta)\ge h_n \quad(\Delta\subset\Lambda_n,\quad\operatorname{rank}\Delta=n-1).$$
If $u\in\Lambda_n^*$ is primitive, the projection of $\Lambda_n$ onto $\mathbb Ru$ has spacing $1/\|u\|_2$. Hence $\det(\Lambda_n\cap u^\perp)=\|u\|_2$, and $\lambda_1(\Lambda_n^*)\ge h_n$. Banaszczyk's Euclidean transference bound~\cite{Banaszczyk1993}, $\mu(B_2^n,\Lambda_n)\lambda_1(\Lambda_n^*)\le n/2$, now yields
$$\frac{\kappa(B_2^n,\Lambda_n)}{\mu(B_2^n,\Lambda_n)} \ge\frac{h_n^2}{2\sqrt n}\ge c\sqrt n.$$

For the upper bound, let
$$V=\{x:\|x\|_2\le\|x-v\|_2\text{ for every }v\in\Lambda\}$$
be the Voronoi cell of $\Lambda$, and put $\mu=\mu(B_2^n,\Lambda)$. For each pair of opposite facets, choose its relevant lattice vector $v$, so that $F_v=V\cap(V+v)$ and its outer unit normal is $v/\|v\|_2$. Since $V=-V$, the map $x\mapsto v-x$ preserves $F_v$. Its centroid is therefore $v/2$, and $V\subset\mu B_2^n$ gives $\|v\|_2\le2\mu$. Applying the divergence theorem to the vector fields $x_i e_j$ and grouping opposite facets gives
\begin{equation}\label{eq:short-voronoi-identity}
    I_n=\sum_v a_vvv^{\mathsf T},\qquad a_v=\frac{\vol_{n-1}(F_v)}{\vol_n(V)\|v\|_2}>0.
\end{equation}

Put $w_v=v/\|v\|_2$ and $p_v=a_v\|v\|_2^2/n$. The vectors $w_v$ are unit vectors, and~\eqref{eq:short-voronoi-identity}, together with its trace, gives
$$\sum_vp_v=1,\qquad \sum_vp_vw_vw_v^{\mathsf T}=I_n/n.$$
The basis-selection result of Harvey and Olver~\cite{HarveyOlver2014} (\href{https://arxiv.org/pdf/1307.2274v1\#page=9}{arXiv v1, Theorem~4.3}) gives independent lattice vectors $b_1,\ldots,b_n$ among the relevant vectors such that
$$\sum_{j=1}^n\frac{b_jb_j^{\mathsf T}}{\|b_j\|_2^2}\preceq C_0I_n,$$
where $C_0$ is absolute and $\preceq$ denotes the positive semidefinite order. Since $\|b_j\|_2\le2\mu$, their column matrix $B$ satisfies
$$BB^{\mathsf T} =\sum_{j=1}^n b_jb_j^{\mathsf T} \preceq4C_0\mu^2I_n.$$
Every Kuhn simplex of this frame lies in its parallelotope $P_B=\sum_j[0,b_j]$. Its center is $c_B=\frac12\sum_jb_j$, and
$$\max_{x\in P_B}\|x-c_B\|_2 =\frac12\max_{\varepsilon\in\{-1,1\}^n}\|B\varepsilon\|_2 \le\frac{\sqrt n}{2}\|B\|_{\mathrm{op}} \le\sqrt{C_0n}\,\mu.$$
Thus the whole parallelotope, and therefore every permutation simplex, fits in one translate of the stated ball. This proves the upper bound. The case $n=1$ is immediate.
\end{proof}

The difference between $\kappa$ and $\mu$ also affects estimates in terms of projection volumes. Define the Kannan--Lov\'asz parameter by
\begin{equation}\label{eq:kl-volume-parameter}
    \mu_{\mathrm{KL}}(K,\Lambda)= \sup_{\substack{1\le d\le n\\L\in\latG(\Lambda,n-d)}} \left(\frac{\det(\Lambda|L^\perp)} {\vol_d(K|L^\perp)}\right)^{1/d}.
\end{equation}
Every lattice covering has density at least one. Applying this fact in each lattice quotient gives $\mu_{\mathrm{KL}}\le\mu$. The lattices in Theorem~\ref{thm:euclidean-separation} therefore satisfy
$$\frac{\kappa(B_2^n,\Lambda_n)} {\mu_{\mathrm{KL}}(B_2^n,\Lambda_n)} \ge c\sqrt n.$$
In particular, for any fixed $a\ge0$, there is no constant $C_a$ independent of dimension such that $\kappa(K,\Lambda)\le C_a\log^a(2n)\mu_{\mathrm{KL}}(K,\Lambda)$, even when $K$ is restricted to Euclidean balls. For the covering radius, Reis and Rothvoss \cite[Theorem~3]{ReisRothvoss2026} prove
$$\mu(K,\Lambda)\le C\log^3(2n)\mu_{\mathrm{KL}}(K,\Lambda).$$
In the next section, we allow more periodic triangulations and obtain minima equal to the covering minima on ellipsoids.

\section{Allowing periodic triangulations}\label{sec:extension}

We relax the common-frame condition by allowing arbitrary periodic lattice triangulations. The resulting minima agree with the covering minima on ellipsoids. We also allow projections of triangulations in higher dimensions.

Let $V$ have any finite dimension $N\ge n$, let $\widetilde\Lambda\subset V$ be a full-rank lattice, and let $\pi:V\to\R^n$ be a linear surjection with $\pi\widetilde\Lambda=\Lambda$. Choose a locally finite geometric triangulation $\mathcal T$ whose simplices cover $V$ and meet in common faces, with vertices in $\widetilde\Lambda$, periodic under a full-rank sublattice $\Gamma\subset\widetilde\Lambda$, and with finitely many maximal simplex orbits. Writing $\mathcal T^{(N)}$ for its $N$-dimensional simplices, define
\begin{equation}\label{eq:extended-definition}
    \hk(K,\Lambda)=\inf_{V,\widetilde\Lambda,\Gamma,\pi,\mathcal T} \max_{[\Delta]\in\mathcal T^{(N)}/\Gamma}r_K(\pi\Delta),
\end{equation}
and define $\hk_i$ by the same quotient formula as \eqref{eq:kappa-minima}. The projected cells may overlap; the complete triangulation in the original space is part of the data. If $\rho$ is a lattice quotient, the map $\rho\circ\pi$ sends $\widetilde\Lambda$ onto $\rho\Lambda$. Thus the same lattice $\widetilde\Lambda$ and triangulation $\mathcal T$, with projection $\rho\circ\pi$, can be used for $(\rho K,\rho\Lambda)$.

\begin{proposition}\label{prop:extension}
For every convex body and lattice,
$$\mu_i(K,\Lambda)\le\hk_i(K,\Lambda)\le\kappa_i(K,\Lambda).$$
For every ellipsoid $E\subset\R^n$,
$$\hk_i(E,\Lambda)=\mu_i(E,\Lambda)\qquad(1\le i\le n).$$
\end{proposition}
\begin{proof}
Kuhn triangulations with the identity projection give the upper bound. For the lower bound, choose a linear right inverse $\iota:\R^n\to V$ of $\pi$ and identify $V$ with $\R^n\oplus\ker\pi$ by $(x,y)\mapsto\iota(x)+y$. Let $B_{\ker\pi}$ be the unit ball for a fixed Euclidean norm on $\ker\pi$, and suppose all projected simplices fit in translates of $tK$. Finitely many simplex orbits give a uniform diameter bound, so for some fixed $D$ each source simplex satisfies
$$\Delta\subset(x_\Delta+tK)\times(y_\Delta+DB_{\ker\pi})$$
for suitable $x_\Delta\in\R^n$ and $y_\Delta\in\ker\pi$. Balitskiy's lemma and projection give
$$\bigl(tK\times DB_{\ker\pi}\bigr)+\widetilde\Lambda=V\quad\Longrightarrow\quad tK+\Lambda=\R^n.$$
Thus $\mu(K,\Lambda)\le t$ for every choice in~\eqref{eq:extended-definition}, and hence $\mu(K,\Lambda)\le\hk(K,\Lambda)$. Applying this inequality and the ordinary upper bound in every lattice quotient gives $\mu_i\le\hk_i\le\kappa_i$.

For the Euclidean unit ball, the Delone cells are convex hulls of lattice points on empty spheres and form a periodic decomposition. Each lies in a ball of radius at most $\mu(B_2^n,\Lambda)$. By Dutour Sikiri\'c and Garber~\cite[Lemma~8]{DutourGarber2020}, this decomposition admits a lattice-periodic triangulation that refines every cell, agrees on common faces, and uses no new vertices. The simplices of this triangulation remain in the original balls and have finitely many orbits: after translating one vertex to zero, every other vertex belongs to the finite set $\Lambda\cap2\mu(B_2^n,\Lambda)B_2^n$. Taking this triangulation with the identity projection proves $\hk(B_2^n,\Lambda)\le\mu(B_2^n,\Lambda)$.

For an invertible linear map $A$, replace $\pi$ by $A\pi$. A translation of $K$ only changes the fitting translations. Hence
$$\hk(AK+x,A\Lambda)=\hk(K,\Lambda)\qquad(x\in\R^n).$$
Every lattice quotient of an ellipsoid is an ellipsoid, so the ordinary equality gives the indexed one.
\end{proof}

For ellipsoids, the equality in Proposition~\ref{prop:extension} retains the known projection-volume estimate \cite[Theorem~73]{ReisRothvoss2026}:
$$\hk(E,\Lambda)=\mu(E,\Lambda) \le C\sqrt{\log(2n)}\,\mu_{\mathrm{KL}}(E,\Lambda).$$
For general convex bodies, this triangulation argument must control all projected simplices at once, while keeping a complete triangulation with common faces in the original space. It remains open whether the following bounds hold.

\begin{question}\label{question:extended-kl}
Is there an absolute constant $C$ such that, for every convex body $K\subset\R^n$ and full-rank lattice $\Lambda$,
$$\hk(K,\Lambda)\le C\log^3(2n)\mu_{\mathrm{KL}}(K,\Lambda)?$$
Can the factor be reduced to $C\log(2n)$?
\end{question}

The first estimate would extend the Reis--Rothvoss bound to this stronger geometric requirement. Since $\mu\le\hk$, the second would also imply the subspace flatness conjecture.

\end{document}